\documentclass{amsart}
\usepackage{amscd,amssymb,amsmath,amsthm}
\usepackage[dvips]{graphicx}
\usepackage[all]{xy}
\theoremstyle{plain}
\newtheorem{proposition}{Proposition}
\newtheorem{theorem}[proposition]{Theorem}

\newtheorem*{proposition*}{Proposition}
\newtheorem*{theorem*}{Theorem}
\newtheorem*{corollary*}{Corollary}
\newtheorem*{lemma*}{Lemma}
\newtheorem*{remark*}{Remark}
\newtheorem*{example*}{Example}
\newtheorem*{definition*}{Definition}

\newcommand{\Q}{\mathbb{Q}}
\newcommand{\R}{\mathbb{R}}

\begin{document}

\title{Obstruction classes and Open Gromov-Witten invariants}

\author{Vito Iacovino}

\address{}

\email{vito.iacovino@gmail.com}

\date{version: \today}


\begin{abstract}
We define genus zero open Gromov-Witten invariants for Calabi-Yau three folds and relatively spin Lagrangian submanifold of Maslov index zero.   
\end{abstract}

\maketitle

\section{Introduction}

Consider a Calabi-Yau 3-fold $X$ and a Maslov index zero Lagrangian submanifold $L$ oriented and relatively spin. 
The open Gromov-Witten potential associated to the pair $(X,L)$ is a solution of a master equation (see \cite{OGW2}). 
In the genus zero case, it follows that the value of the potential on its critical points lead to invariants (see \cite{OGW2}).     
In this paper we give an alternative direct definition of these invariants that does not use the potential. 
We extend the method used in \cite{OGW1}, where the particul case with $L$ of rational homology of a sphere is considered.

The invariants associated to the Gromov-Witten potential are defined as real numbers and are associated to a given symplectic area. 
With the method of this paper the invariants are manifestly rational numbers and 
we can also associate invariants to relative homology classes of the pair $(X,L)$. 

If we try to extend the argument of \cite{OGW1}, that uses generalized linking numbers, we are naturally conducted to consider obstruction classes. 
We find that open Gromov-Witten invariants have obstructions and ambiguities analogous to the one found in \cite{FO3} for Lagrangian Floer Homology.
In fact this is more than just an analogy. 
The data that we need in this paper to define open Gromov-Witten invariants are strictly related to the bounding chains of \cite{FO3}.
We use this in \cite{floer} to provide a new construction of Lagrangian Floer Homology.
This applies in any dimension even dropping the Calabi-Yau assumption.
Our approach is thereby simpler compared to the approach of \cite{FO3} since we do not use any chain level structure.

\section{Obstruction and Ambiguity of Open Gromov-Witten Invariants}

Let $X $ be a Calabi-Yau and $L$ be relative spin Lagrangian submanifold. 
We consider the problem of counting disks in some relative homology class $A \in H_2(X,L)$. 
More precisely we want to count multi-disks in the sense of \cite{OGW1}. 

As done in \cite{OGW1} in the special case that $L$ has rational homology of a sphere, 
we consider the problem recursively starting with disks of minimal area.

For the relative homology class $A_0$ of minimal area the number of disks $\mathcal{F}_0(A_0)$ is well defined 
since in this case the relevant moduli space has no boundary.
Define 
\begin{equation} \label{obs1}
o(A_0) =  W_{T_{0,1}(A_0)} \in C_1(L)
\end{equation}
where $W_{T_{0,1}(A_0)}$ is the singular chain associated to the tree with one vertex in the the homology class $A_0$ and one external edge (\cite{OGW1}).
In the next step we will see that this is the first obstruction to define open Gromov-Witten invariants. 
Observe that $o(A_0)$ is defined up to homotopy in the sense of \cite{OGW1}. 
Two representative of $o(A_0)$ are connected by an homotopy uniquely determined up to equivalence.

Let us now consider the next step $F_0(2 A_0)$. 
The moduli space $\overline{\mathcal{M}}_0(2 A_0)$ has a boundary therefore we need to consider the moduli space of multi-disks. 
Following \cite{OGW1} the associated system of chains can be written as
\begin{equation} \label{system1}
\{ W_{T_0} , W_{T_1}  \}.   
\end{equation}
where $T_0$ is the tree with one vertex in the homology class $2 A_0$ and 
$T_1$ is the tree given by two vertices in the relative class $A_0$ connected by an edge. 
In order to associate a number to (\ref{system1}) we need to contract the singular chain $ W_{T_1} \in C_2(L^2)$. 
By the gluing property there exists an homotopy (uniquely determined up to equivalence) of $W_{T_1} $ with 
$  W_{T_{0,1}(A_0)} \times  W_{T_{0,1}(A_0)}   $. 

 

Now we need to make the assumption that the first obstruction class (\ref{obs1}) vanishes as element of $H_1(L, \Q)$. 
In this case there exists $b (A_0) \in C_2(L)$ such that
$$ \partial b (A_0) = o(A_0) .$$
We consider another $b' (A_0)$ solution of $ \partial b' (A_0) = o(A_0) $ equivalent to $b(A_0)$ if $b' (A_0) - b(A_0)$ is zero in $H_2(L, \Q)$.
We can use $b (A_0)$ to define an homotopy of $ W_{T_{0,1}(A_0)} $ with the zero singular chain. 
Applying the procedure of \cite{OGW1} to the system of chains $(\ref{system1})$ we get a rational number $F_0(2 A_0)$. 
This number will depend on the equivalence class of $b (A_0)$.

We iterate the above process.
Let $A \in H_2(X,L)$, and suppose that we defined $o(A')$ and $b(A')$ with $ \partial b (A') = o(A')$
for all $A'$ of symplectic area strictly less than the symplectic area of $A$. 
Using the procedure of \cite{OGW1} we can find an homotopy between $W_{\mathcal{T}(A)}$ and a system of chains 
$W_{\mathcal{T}(A)}^0$ such that 
$$ W_T^0 =0  $$
if $T$ has at least an internal edge.
The equivalence class of this homotopy depends on the equivalence class of the $b$.
Define the open Gromov-Witten invariant  
\begin{equation} \label{invariant}
F_0(A) =  W^0_{T_0(A)} \in \Q
\end{equation}
where $T_0(A)$ is the tree in the the homology class $A$ with exactly one vertex and no edges.
Define also the obstruction class as 
\begin{equation} \label{obs}
o(A) =  W^0_{T_{0,1}(A)} \in C_1(L)
\end{equation}
where $T_{0,1}(A)$ is the the tree in the the homology class $A$ with exactly one vertex, no internal edges and one external edge.

Observe that in $ H_1(L)$
$$ [o(A)] = F_0(A) [\partial A] .$$
Therefore the condition that the obstruction class vanishes is equivalent to the vanishing of open Gromov-Witten 
in the case that the boundary wind a not trivial cycle.

\begin{theorem}
For each $A \in H_2(X,L)$ there exist 
$$   o(A)  \in C_1(L, \Q)$$
$$   b(A)  \in C_2(L, \Q)$$
such that
\begin{enumerate}
\item $o(A)$ is defined when $o(A')$ and $b(A')$ are defined for all $A' \in H_2(X,L)$ such that $ \omega(A ) < \omega(A ') $.
In this case also the open Gromov-Witten invariant $F_0(A)$ counting disks in the class $A$ is defined. 
\item $b(A)$ is defined if $o(A)$ is defined and $[o(A)]=0$ in $H_1(L, \Q)$. $b(A)$ is a singular chain such that $ \partial B (A) =  o (A) .$ 
\end{enumerate}
\end{theorem}

\subsection{Invariants associated to symplectic area}

The vanishing of the obstruction classes associated to relative homology classes is a restrictive condition.
We now consider the problem of counting disks of some symplectic area. 
This are less refined invariants, but are defined with weaker conditions.

In this subsection we consider decorated trees with vertices associated with a symplectic area instead of a relative homology class.
 
The method of \cite{OGW1} applies as well in this context after the following observation. 
In order to define $F_0(E)$ and $o(E)$ we are interested only to trees with at most one external edge.
In this case all the relevant trees have at least one vertex of valence one, unless the tree is the trivial tree with one vertex and no edges. 
From this it follows that in the recursive argument above we can restrict to contract only components associated to vertices of valence one.

We have then (compare Theorem $3.6.18$ of \cite{FO3})
\begin{theorem}
There is a sequence of $o(E_k)$, $b(E_k)$ for $E_k \in \R$ with 
$$   o(E_k)  \in C_1(L, \Q)$$
$$   b(E_k)  \in C_2(L, \Q)$$
such that
\begin{enumerate}
 \item $o(E_k)$ is defined when $o(E_{k'})$ and $b(E_{k'})$ are defined for all $k' < k$.
In this case also the open Gromov-Witten invariant $F_0(E_k)$ counting multi-disks of energy $E_k$ is defined. 
$F_0(E_k)$ depends on the choices of the equivalent class of $b(E_{k'})$.
\item $b(E_k)$ is defined if $o(E_k)$ is defined and $[o(E_k)]=0$ in $H_1(L, \Q)$. $b(E_k)$ is a singular chain such that $ \partial b (E_k) =  o (E_k) .$ 
\end{enumerate}
\end{theorem}

\subsection{Anti-simplectic involutions}
We now consider the case that $L$ arises as the fixed points of an anti-symplectic involution $\phi$.
In this case we show that all the obstruction classes vanishes.

For each decorated tree $T$ and for each vertex $v \in V(T)$ there is a natural involution on $ \overline{\mathcal{M}}_T $ inducted by the action of $\phi$
on the component associated to $v$. 
In analogy with \cite{So}, consider perturbations invariants by all of these involutions. 
The associate system of chains $W_{\mathcal{T}}$ has the property that 
$$ W_T = 0$$
if $T$ has at least a vertex of valence one.
Therefore $ W_T = 0$ if $T$ has exactly one external edge. 
It follows that all obstruction classes vanish and we can set $b(E)=0$ for all $E$.
The associate open Gromov-Witten invariant coincide with the one of \cite{So}.


\end{document}